\documentclass[11pt]{article}
\usepackage[colorlinks=true,linkcolor=blue,citecolor = blue]{hyperref}
\usepackage {verbatim}
\usepackage[commandnameprefix=always]{changes}
\usepackage{mathrsfs}
\usepackage{amsfonts}
\usepackage{amssymb}
\usepackage{amsthm}
\usepackage{amsmath}
\usepackage{graphicx}
\usepackage{etoolbox,xstring,mfirstuc,textcase}
\usepackage{empheq}
\usepackage{indentfirst}
\usepackage{cite} 
\usepackage{cases}
\usepackage{graphics}
\usepackage{xcolor}  
\usepackage{amsmath,bm}
\usepackage{changes}
\usepackage{lineno}
\usepackage{cases}
\usepackage{hyperref}
\newtheorem{theorem}{\textbf{Theorem}}[section]
\newtheorem{lemma}{\textbf{Lemma}}[section]
\newtheorem{proposition}{\textbf{Proposition}}[section]
\newtheorem{corollary}{\textbf{Corollary}}[section]
\newtheorem{remark}{\textbf{Remark}}[section]
\newtheorem{definition}{\textbf{Definition}}[section]
\allowdisplaybreaks[4] 
\def\be{\begin{equation}}
	\def\ee{\end{equation}}
\def\bt{\begin{theorem}}
	\def\et{\end{theorem}}
\def\bl{\begin{lemma}}
	\def\el{\end{lemma}}
\def\br{\begin{remark}}
	\def\er{\end{remark}}
\def\bp{\begin{proposition}}
	\def\ep{\end{proposition}}
\def\bc{\begin{corollary}}
	\def\ec{\end{corollary}}
\def\bd{\begin{definition}}
	\def\ed{\end{definition}}

\begin{document}
	
\title{Estimates on the Laplace Operator in Heat Flows of  Harmonic Maps }
	
\author{
	{Qingtong Wu}
	\footnote{School of Mathematical Sciences, Fudan University, Handan Road 220, Shanghai 200433, China. Email:  \textit{qtwu22@m.fudan.edu.cn}
        }
}

\date{}
	
\maketitle
	

\begin{abstract}
\noindent 
In this paper
we investigate estimates about the Laplace operator in heat flows of harmonic maps, focusing outside the singularities through spherical coordinates. These estimates can be used in the general Ericksen--Leslie system to obtain higher-order estimates. 
We consider the problem subject to the $\mathbb{T}^2$ and $\mathbb{T}^3$ boundary conditions.
  \medskip \\
\medskip\noindent

\end{abstract}
\section{Introduction}
The dynamic continuum theory for Liquid crystal flows was established by Leslie (1968)\cite{ref5} and Ericksen (1961)\cite{ref1}.
$\mathbf{d}:\mathbb{T}^n\rightarrow\mathbb{S}^2$  $(n=2,3)$
is the director field describing the average orientation of a neighboring region
and
satisfies the following boundary condition
\begin{equation}
\mathbf{d}(x+\mathbf{e}_i,t)
=
(-1)^{a_i}\mathbf{d}(x,t)    ,  \; 
\text{for} \;
(x,t)\in \partial Q\times(0,\infty).
\end{equation}
Unit vectors $\mathbf{e}_i \;(i=1, \ldots,n)$ are the canonical basis of $\mathbb{R}^n$ and
$Q=\prod\limits_{i=1}^{n}[0,1]$ is a unit square in $\mathbf{R}^n$.
The more general case $Q=\prod\limits_{i=1}^{n}[0,L_i]$
$(n=2,3)$
with different periods $L_i$ in different directions $\mathbf{e}_i$ can be treated in a same way.
The integers \(a_i\; (i=1, \ldots,n)\)  are constants depending on
\[
\mathbf{d}_0(x+\mathbf{e}_i) = (-1)^{a_i}\mathbf{d}_0(x), \;
x\in \partial Q.
\]
This is because the liquid crystal molecules coincide with their original positions after rotating by an angle of \(k_i\pi\) around the center, where the integers \(k_i\) \((i=1, \ldots,n)\) are constants.

In this paper, we study the estimates on the Laplace operator  
in heat flows of 
harmonic maps,
which can be used in
the global well-posedness of the general Ericksen–Leslie system and heat flows of harmonic maps.
Denote $\mathbf{d}=(d_1,d_2,d_3)^T$.
Lie, Li and Zhang \cite{ref4}
considered the Cauchy problem to 
the two-dimensional incompressible liquid crystal equation and the heat flows of the harmonic maps equation.
They got a similar 
conclusion in Theorem 2.1 
when there exists $\epsilon_0 >0$ $s.t.$ $d_3>\epsilon_0$.

Inspired by \cite{ref3}, 
denote
$\mathbf{n}=(\sin\phi,\cos\phi)^T$
and
$$\mathbf{d}=(\mathbf{n}\sin\theta,\cos\theta)^T=(\sin\phi\sin\theta,\cos\phi\sin\theta,\cos\theta)^T,$$
then
$\mathbf{n}^\perp=(\cos\phi,-\sin\phi)^T$, 
$\mathbf{d}^\perp_1=(\mathbf{n}^\perp,0)^T$
and
$\mathbf{d}^\perp_2=(\mathbf{n}\cos\theta,-\sin\theta)^T$,
where
$\phi=\phi(x_1,x_2)$
and
$\theta=\theta(x_1,x_2)$.
It is obvious that 
$\mathbf{d}^\perp_1\cdot\mathbf{d}^\perp_2
=\mathbf{d}^\perp_1\cdot\mathbf{d}
=\mathbf{d}\cdot\mathbf{d}^\perp_2
=0$.
$$
\begin{aligned}
|\nabla\mathbf{d}|^2
=&\sin^2\theta |\nabla\phi|^2 +  |\nabla\theta|^2,\\
\Delta\mathbf{d}
=& \left( \Delta\phi \sin\theta
+\nabla\phi\cdot\nabla\theta \cos\theta 
\right)  \mathbf{d}^\perp_1
+\Delta\theta \mathbf{d}^\perp_2
+( \sin\theta \nabla\phi \cdot \nabla\mathbf{d}^\perp_1
+\nabla\theta \cdot \nabla\mathbf{d}^\perp_2
).
\end{aligned}
$$
By calculating directly,
$$
\sin\theta \nabla\phi \cdot \nabla\mathbf{d}^\perp_1
+\nabla\theta \cdot \nabla\mathbf{d}^\perp_2
+|\nabla\mathbf{d}|^2 \mathbf{d}
=\nabla\phi \cdot \nabla\theta \cos\theta \mathbf{d}^\perp_1
-\sin\theta\cos\theta |\nabla\phi|^2 \mathbf{d}^\perp_2.
$$
Then we have the following identity in two versions with two different meanings.
\begin{equation}
\Delta\mathbf{d}+|\nabla\mathbf{d}|^2\mathbf{d}
= \left( \Delta\phi \sin\theta
+2\nabla\phi\cdot\nabla\theta \cos\theta 
\right)  \mathbf{d}^\perp_1
+(\Delta\theta- \sin\theta\cos\theta |\nabla\phi|^2 )\mathbf{d}^\perp_2,   
\end{equation}
\begin{equation}
\Delta\mathbf{d}
=-|\nabla\mathbf{d}|^2\mathbf{d}
 +\left( \Delta\phi \sin\theta
+2\nabla\phi\cdot\nabla\theta \cos\theta 
\right)  \mathbf{d}^\perp_1
+(\Delta\theta- \sin\theta\cos\theta |\nabla\phi|^2 )\mathbf{d}^\perp_2.    
\end{equation}

\begin{remark}
(i)
The above derivation considers the polar angle with respect to the z-axis,
i.e. $(\theta,\phi)=(\theta^z,\phi^z)$.
We can denote 
$$
\begin{aligned}
\mathbf{d}=&(\sin\phi^z\sin\theta^z,\cos\phi^z\sin\theta^z,\cos\theta^z)^T\\
=&(\cos\theta^x,\sin\phi^x\sin\theta^x,\cos\phi^x\sin\theta^x)^T\\
=&
(\cos\phi^y\sin\theta^y,\cos\theta^y,\sin\phi^y\sin\theta^y)^T.
\end{aligned}
$$
(ii)
The basic energy law  of the general  Ericksen–Leslie system \cite{ref2}\cite{ref7}
can be written as follows,
\begin{equation}
\begin{aligned}
&\frac{1}{2} \frac{\mathrm{d}}{\mathrm{d} t} 
\int_Q
|\mathbf{v}|^2
+\frac{1-\gamma}{\operatorname{Re}}
(\sin^2\theta |\nabla\phi|^2 +  |\nabla\theta|^2
 ) 
\mathrm{~d} x\\
=&  
-\int_Q
\frac{\gamma}{\operatorname{Re}}
|\nabla \mathbf{v}|^2
+\frac{1-\gamma}{\operatorname{Re}}[\beta_1(\mathbf{d}\otimes\mathbf{d}: \mathbf{D})^2
+\beta_2 \mathbf{D}: \mathbf{D} 
+\beta_3|\mathbf{D}\mathbf{d}|^2
+\mu_1|\mathbf{h}|^2
-\mu_1(\mathbf{h} \cdot \mathbf{d})^2] 
\mathrm{d} x .
\end{aligned}
\end{equation}	
where 
$\mathbf{v}:\mathbb{T}^n \rightarrow \mathbb{R}^n$ 
is the velocity of fluids, 
$\mathbf{h}=\Delta\mathbf{d}$.
Re is the Reynolds number, 
 $\gamma \in(0,1)$ and $0<\mu_1$. 
D represents the rate of strain tensor and $\boldsymbol{\Omega}$ is the skew-symmetric part of the strain rate
$$
\mathbf{D}=\frac{1}{2}\left[\nabla \mathbf{v}+(\nabla \mathbf{v})^T\right]=\frac{1}{2}\left(\frac{\partial \mathbf{v}_i}{ \partial x_j}+\frac{\partial \mathbf{v}_j}{\partial x_i}\right), \;
\boldsymbol{\Omega}=\frac{1}{2}\left[\nabla \mathbf{v}-(\nabla \mathbf{v})^T\right]=\frac{1}{2}\left(\frac{\partial \mathbf{v}_i}{\partial x_j}-\frac{\partial \mathbf{v}_j}{\partial x_i}\right) .
$$	

 $(\theta,\phi)$ can be replaced by 
$(\theta^z,\phi^z)$,
$(\theta^y,\phi^y)$
and
$(\theta^x,\phi^x)$.
Then 
$
\|\nabla\mathbf{d}\|^2=\sin^2\theta |\nabla\phi|^2 +  |\nabla\theta|^2    
$ 
is bounded if the initial conditions $(\mathbf{v}_0,\mathbf{d}_0)$ is good enough.
Without special instructions, we take 
$(\theta,\phi)$ to be 
$(\theta^z,\phi^z)$.
\end{remark}

We need the following  Gagliardo–Nirenberg inequality which has the same solution in $\mathbb{T}^2$ and $\mathbb{T}^3$  boundary conditions.
\begin{lemma}[Gagliardo–Nirenberg Inequality \cite{ref6}]
Let \( \Omega \subset \mathbb{R}^n \) be a bounded domain with sufficient smoothness, and let \( u \in H^m(\Omega) \cap L^p(\Omega) \) for some \( m \in \mathbb{N} \) and \( 1 \leq p \leq \infty \). Then there exists a constant \( C \) depending only on \( \Omega \), \( m \), \( p \), and \( n \), such that
\begin{equation}
 \| u \|_{L^q(\Omega)} \leq C \| u \|_{H^m(\Omega)}^\theta \| u \|_{L^p(\Omega)}^{1-\theta}, 
\end{equation}
where \( \theta = \frac{n \left( \frac{1}{p} - \frac{1}{q} \right) + m}{n + m} \) and \( q \) is the Sobolev conjugate exponent of \( p \), i.e. \( \frac{1}{q} = \frac{1}{p} - \frac{m}{n} \).
\end{lemma}

\section{Estimates}
\setcounter{equation}{0}
\begin{theorem}
Assuming that $\mathbf{d}:\mathbb{T}^2\rightarrow\mathbb{S}^2$ and
$\mathbf{d}\in C^{\infty}$
, 
we have the following estimate,
\begin{equation}
    \|\nabla^2\mathbf{d}\|
\leq
    C
    \left(
    \| \Delta\mathbf{d}+|\nabla\mathbf{d}|^2\mathbf{d} \|
    +1
    \right)
\end{equation}
where $C$ is a constant depending on $\|\nabla\mathbf{d}\|$.
\end{theorem}

\paragraph{Proof.}
(i) First, we have
$
 \|\nabla^2\mathbf{d}\|^2
=\int_Q \Delta\mathbf{d}\cdot \Delta  \mathbf{d} 
\mathrm{~d}x
$ 
and
\begin{equation}
  \|  \Delta\mathbf{d}  \|^2
=
 \| \Delta\mathbf{d}+|\nabla\mathbf{d}|^2\mathbf{d} \|^2
+\|\nabla\mathbf{d}\|^4_{L^4},
\end{equation}
where
$$ \|\nabla\mathbf{d}\|^4_{L^4}
=
\|  \sin^2\theta |\nabla\phi|^2 +  |\nabla\theta|^2   \|^2
\leq
2
\left(
\|\sin\theta \nabla\phi\|^4_{L^4}
+
\|\nabla\theta\|^4_{L^4}
\right).
$$

Let
$$
\begin{aligned}
  \nabla\cdot\left(\sin\theta \nabla\phi \right)
+\nabla\phi\cdot\nabla\theta \cos\theta 
&= g_1,
\\
\Delta\theta- \sin\theta\cos\theta |\nabla\phi|^2 
&=g_2.
\end{aligned}
$$
It is obvious that $\|g_1\|^2+\|g_2\|^2= \left\| \Delta\mathbf{d}+|\nabla\mathbf{d}|^2\mathbf{d} \right\|^2$,
and we have
\begin{equation}
\begin{pmatrix} 
\nabla \cdot ( \sin\theta\nabla\phi) \\ 
\nabla \cdot  \nabla \theta
\end{pmatrix} 
+ 
\begin{pmatrix} 
\nabla \theta  \cdot 
\left(\cos\theta \nabla\phi\right)
\\ 
-\sin \theta \nabla \phi \cdot 
\left(\cos\theta \nabla\phi\right)
\end{pmatrix} 
= 
\begin{pmatrix} 
g_1 \\ g_2 
\end{pmatrix}.
\end{equation}
Let $\mathbf{u}=(\sin\theta\nabla\phi ,\nabla \theta)^T$,
then we have
\begin{equation}
    \|\nabla\cdot(\sin\theta\nabla\phi)\|^2
    +
    \| \nabla\cdot(\nabla \theta) \|^2
    \leq
    \|\mathbf{u}\|_{L^4}^2 \| \cos\theta \nabla\phi \|_{L^4}^2
    +\|g_1\|+\|g_2\|.
\end{equation}
By the following Gagliardo–Nirenberg inequality\cite{ref6}, 
for $f:\mathbb{T}^2 \rightarrow \mathbb{R}^2$,
$$
\| f\|^4_{L^4}
\leq
C\|f\|^2 
\left(
\|\nabla f\|^2+\|f\|^2
\right),
$$
from which we can only obtain (2.1)
under small initial conditions and 
$\mathbf{d}$ which can not reach z-axis.

For $\mathbf{u}_1:\mathbb{R}^2\rightarrow\mathbb{R}^2$,
$\mathbf{u}_2:\mathbb{R}^2\rightarrow\mathbb{R}^2$,
$\mathbf{b}=\cos\theta\nabla\phi$,
$\mathbf{u}_1=\sin\theta\nabla\phi$,
$\mathbf{u}_2=\nabla\theta$.
Denote
\begin{equation}
\mathbf{u}:=
\begin{pmatrix}
&\mathbf{u}_1 \\
&\mathbf{u}_2
\end{pmatrix}
,\quad
\mathbf{A} \mathbf{u}:=
\begin{pmatrix}
&\nabla\cdot\mathbf{u}_1 \\
&\nabla\cdot\mathbf{u}_2
\end{pmatrix}
,\quad
\mathbf{B} \mathbf{u}:=
\begin{pmatrix}
&\mathbf{b}\cdot\mathbf{u}_2 \\
&-\mathbf{b}\cdot\mathbf{u}_1
\end{pmatrix},
\end{equation}
where $\mathbf{A} := \prod\limits_{j=1}^{4} C^{\infty} (Q) \rightarrow \prod\limits_{j=1}^{2}C^{\infty} (Q)  $,
$\mathbf{B}:=  \prod\limits_{j=1}^{4} C^{\infty} (Q) \rightarrow \prod\limits_{j=1}^{2}C^{\infty} (Q) $.
Obviously $\mathbf{A}$,
$\mathbf{B}$
are linear operators.
Then we just need to prove that there exists a constant $M>0$ $s.t.$ 
\begin{equation}
\|\mathbf{A} \mathbf{u}\|^2
\leq
M \left(\|(\mathbf{A}+\mathbf{B}) \mathbf{u}\|^2  +1  \right).
\end{equation}

Let 
$ D_1:= \left((M-1)\mathbf{A}\mathbf{u},\mathbf{A}\mathbf{u} \right)
+
\left(M (2\mathbf{A}+\mathbf{B}) \mathbf{u},
\mathbf{B}\mathbf{u}\right) 
=
M\|\Delta\mathbf{d}+|\nabla\mathbf{d}|\|^2 
-
\| \mathbf{A}\mathbf{u} \|^2
$,
then we have
\begin{equation}
\begin{aligned}
D_1=
&\left((M-1)\mathbf{A}\mathbf{u},\mathbf{A}\mathbf{u} \right)
+
\left(M (2\mathbf{A}+\mathbf{B}) \mathbf{u},
\mathbf{B}\mathbf{u}\right)\\
=
&
(M-1)
\int_Q
|\nabla\cdot\mathbf{u}_1|^2
+|\nabla\cdot\mathbf{u}_2|^2
\mathrm{~d} x
+
2M
\int_Q
(\nabla\cdot\mathbf{u}_1)(\mathbf{b}\cdot\mathbf{u}_2)
-(\nabla\cdot\mathbf{u}_2)(\mathbf{b}\cdot\mathbf{u}_1)
\mathrm{~d} x\\
&+
M
\int_Q
|\mathbf{b}\cdot\mathbf{u}_1|^2
+
|\mathbf{b}\cdot\mathbf{u}_2|^2
\mathbf{~d} x.
\end{aligned}
\end{equation}

The following term is essential,
\begin{equation}
\begin{aligned}
\int_Q
(\nabla\cdot\mathbf{u}_1)(\mathbf{b}\cdot\mathbf{u}_2)
-(\nabla\cdot\mathbf{u}_2)(\mathbf{b}\cdot\mathbf{u}_1)
\mathrm{~d} x.
\end{aligned}
\end{equation}

Let
$\Theta=\theta+\frac{\pi}{2}$.
We do not need to consider its geometric meaning, but we should choose a right-handed coordinate system to avoid the complications that directional changes due to axial symmetry might introduce in the operations of differential operators.
Recalculate the aforementioned process once more.
Then we have
$$
\begin{aligned}
\widetilde{\mathbf{b}}&=\cos\Theta\nabla\phi =-\sin\theta\nabla\phi,\\
\widetilde{\mathbf{u}}_1&=\sin\Theta\nabla\phi=\cos\theta\nabla\phi,\\
\widetilde{\mathbf{u}}_2&=\nabla\Theta=\nabla\theta,\\
\widetilde{\mathbf{u}} &=
\begin{pmatrix}
\widetilde{\mathbf{u}}_1\\
\widetilde{\mathbf{u}}_2
\end{pmatrix}.
\end{aligned}
$$
Therefore,
\begin{equation}
\begin{aligned}
D_2:=& \left((M-1)\mathbf{A}\widetilde{\mathbf{u}},\mathbf{A}\widetilde{\mathbf{u}} \right)
+
\left(M (2\mathbf{A}+\mathbf{B}) \widetilde{\mathbf{u}},
\mathbf{B}\widetilde{\mathbf{u}}\right) \\
=&
M\|\Delta\mathbf{d}+|\nabla\mathbf{d}|\|^2 
-
\| \mathbf{A}\widetilde{\mathbf{u}}\|^2\\
=
&
(M-1)
\int_Q
|\nabla\cdot\widetilde{\mathbf{u}}|^2
+|\nabla\cdot\widetilde{\mathbf{u}}|^2
\mathrm{~d} x
+
2M
\int_Q
(\nabla\cdot\widetilde{\mathbf{u}}_1)(\widetilde{\mathbf{b}}\cdot\widetilde{\mathbf{u}}_2)
-(\nabla\cdot\widetilde{\mathbf{u}}_2)(\widetilde{\mathbf{b}}\cdot\widetilde{\mathbf{u}}_1)
\mathrm{~d} x\\
&+
M
\int_Q
|\widetilde{\mathbf{b}}\cdot\widetilde{\mathbf{u}}_1|^2
+
|\widetilde{\mathbf{b}}\cdot\widetilde{\mathbf{u}}_2|^2
\mathbf{~d} x,
\end{aligned}
\end{equation}
and
\begin{equation}
(\nabla\cdot\widetilde{\mathbf{u}}_2)(\widetilde{\mathbf{b}}\cdot\widetilde{\mathbf{u}}_1)
+
(\nabla\cdot\mathbf{u}_2)(\widetilde{\mathbf{b}}\cdot\mathbf{u}_1)
=0.
\end{equation}

It is obvious that
$$
\begin{aligned}
\nabla\cdot\mathbf{u}_1
=&
\sin\theta \Delta\phi
+
\cos\theta
\nabla\theta\cdot\nabla\phi,\\
\nabla\cdot\mathbf{b}
=&
\cos\theta\nabla\phi
-\sin\theta\nabla\phi\cdot\nabla\theta.
\end{aligned}
$$
Then we have
$$
\begin{aligned}
\sin\theta(\nabla\cdot \mathbf{u}_1)  
+
\cos\theta (\nabla\cdot \mathbf{b})
=&
\Delta\phi,\\
\cos\theta(\nabla\cdot \mathbf{u}_1)  
-
\sin\theta (\nabla\cdot \mathbf{b})
=&
\nabla\phi \cdot \nabla\theta.
\end{aligned}
$$

Finally,
\begin{equation}
\begin{aligned}
&\int_Q
(\nabla\cdot\mathbf{u}_1)(\mathbf{b}\cdot\mathbf{u}_2)
+
(\nabla\cdot\widetilde{\mathbf{u}}_1)(\widetilde{\mathbf{b}}\cdot\widetilde{\mathbf{u}}_2)
\mathrm{~d} x    \\
=&
\int_Q
(\nabla\cdot\mathbf{u}_1)( \nabla\phi \cdot \nabla\theta) \cos\theta
-
(\nabla\cdot\mathbf{b})( \nabla\phi \cdot \nabla\theta) \sin\theta
\mathrm{~d} x  \\
=&
\int_Q
\left[
\cos\theta(\nabla\cdot \mathbf{u}_1)  
-
\sin\theta (\nabla\cdot \mathbf{b})
\right]^2
\mathrm{~d} x .
\end{aligned}
\end{equation}
Therefore,
\begin{equation}
\|\mathbf{A}\mathbf{u}\|^2
+
\|\mathbf{A}\widetilde{\mathbf{u}}\|^2
\leq
2 \|\Delta\mathbf{d}+|\nabla\mathbf{d}|\|^2,
\end{equation}
and
\begin{equation}
\|\Delta\theta\|^2
\leq
\|\Delta\mathbf{d}+|\nabla\mathbf{d}|\|^2 .
\end{equation}

By Remark 1.1,
we can use the same method to control $\|\Delta\theta^x\|^2$ and $\|\Delta\theta^y\|^2$.

(ii) We prove $|\sin\theta^z\nabla\phi^z|\leq C (|\nabla\theta^x|+|\nabla\theta^y|)$,
which can be deduced from the following inequality.
For $\forall\delta \in \mathbb{R}^2$ small enough
\begin{equation}
  \left|\sin\theta^z (x)  \frac{\phi^z (x+\delta)-\phi^z (x) }{|\delta|}
  \right|
  \leq C 
  \left(\left|\frac{\theta^x (x+\delta)-\theta^x (x) }{|\delta|}\right|
  +
  \left|\frac{\theta^y (x+\delta)-\theta^y (x) }{|\delta|}
  \right|\right)  .
\end{equation}

If $\sin\theta^z (x)=\sin\theta^z (x+\delta)$, 
since the ratio of the length of the minor arc on the circumference to its corresponding chord does not exceed 
$\frac{\pi}{2}$,
Consider the inequality in the x-0-y plane.
\begin{equation}
\begin{aligned}
      &\left|\sin\theta^z (x)  \left(\phi^z (x+\delta)-\phi^z (x) \right)
       \right|\\
       \leq &
      \frac{ \pi}{2} 
      \sqrt{ \left( \cos\theta^x (x+\delta)-\cos\theta^x (x) \right)^2+ \left(\cos\theta^y (x+\delta)-\cos\theta^y (x) \right)^2  }\\
      \leq &
      C 
  \left(\left|\theta^x (x+\delta)-\theta^x (x) \right|
  +
  \left|\theta^y (x+\delta)-\theta^y (x) 
  \right|\right)  .
\end{aligned}
\end{equation}

If $\sin\theta^z (x)  \neq \sin\theta^z (x+\delta)$, 
without loss of generality
, we can assume $|\sin\theta^z (x)| > |\sin\theta^z (x+\delta)|$.
Consider the case of $\sin\theta^z (x)$
and the inequality in the $x-o-y$ plane.
$\left|\sin\theta^z (x)  \left(\phi^z (x+\delta)-\phi^z (x) \right)\right|$ is the arc length corresponding to two angles on a circle with radius $|\sin\theta^z (x)|$,
corresponding to two points $A$ and $B$ on the circle.
Denote the distance from point $A$ to the line $OB$ to be $s$.
\begin{equation}
\begin{aligned}
      &\left|\sin\theta^z (x)  \left(\phi^z (x+\delta)-\phi^z (x) \right)
       \right|
       \leq  \frac{ \pi}{2} s
       \\
       \leq &
      \frac{ \pi}{2} 
      \sqrt{ \left( \cos\theta^x (x+\delta)-\cos\theta^x (x) \right)^2+ \left(\cos\theta^y (x+\delta)-\cos\theta^y (x) \right)^2  }\\
      \leq &
      C 
  \left(\left|\theta^x (x+\delta)-\theta^x (x) \right|
  +
  \left|\theta^y (x+\delta)-\theta^y (x) 
  \right|\right)  .  \Box
\end{aligned} 
\end{equation}

\begin{remark}
Under the assumption of Theorem 2.1, if $\mathbf{d}$
is a harmonic map,
then
$$
\Delta\mathbf{d}
+
|\nabla\mathbf{d}|^2\mathbf{d}=0.
$$
Therefore, we can achieve better results regarding (2.1) and (2.17).
\end{remark}

\begin{theorem}
Assuming that $\mathbf{d}:\mathbb{T}^2\rightarrow\mathbb{S}^2$ and
$\mathbf{d}\in C^{\infty}$
, 
we have the following estimate,
    \begin{equation}
    \|\nabla^3\mathbf{d}\|
\leq
    C\left(
    \| \Delta \mathbf{d} + |\nabla\mathbf{d}|^2 \mathbf{d} \|^2
    +
    \| \nabla(\Delta\mathbf{d}+|\nabla\mathbf{d}|^2\mathbf{d}) \|
    +1 \right),
\end{equation}
where $C$ is a constant depending on $\|\nabla\mathbf{d}\|$.
\end{theorem}

\paragraph{Proof.}
we can obtain 
$
\|\nabla^3\mathbf{d}\|^2
=\int_Q  (\nabla\Delta \mathbf{d}) : (\nabla\Delta \mathbf{d})
\mathbf{~d} x
$
and the following identities,
\begin{equation}
\begin{aligned}
&\|\nabla\Delta\mathbf{d}\|^2\\
=&
\left\|
\nabla
\left[
\left( \Delta\phi \sin\theta
+2\nabla\phi\cdot\nabla\theta \cos\theta 
\right)   \mathbf{d}^\perp_1
+(\Delta\theta- \sin\theta\cos\theta |\nabla\phi|^2 )\mathbf{d}^\perp_2
-|\nabla\mathbf{d}|^2\mathbf{d}
\right]
\right\|^2\\
=&
\left\|
\nabla\left( \Delta\phi \sin\theta
+2\nabla\phi\cdot\nabla\theta \cos\theta 
\right)
+
(\Delta\theta- \sin\theta\cos\theta |\nabla\phi|^2 )
\cos\theta \nabla\phi
+|\nabla\mathbf{d}|^2 \sin\theta\nabla\phi
\right\|^2 \\
+&
\left\|
\left[
\nabla(\Delta\theta- \sin\theta\cos\theta |\nabla\phi|^2 )
-\left( \Delta\phi \sin\theta
+2\nabla\phi\cdot\nabla\theta \cos\theta 
\right)\cos\theta\nabla\phi
+|\nabla\mathbf{d}|^2 \nabla\theta
\right]
\otimes \mathbf{d}^\perp_2
\right\|^2\\
+&
\left\|
\left[
(\Delta\theta- \sin\theta\cos\theta |\nabla\phi|^2 )\nabla\theta
+
\left( \Delta\phi \sin\theta
+2\nabla\phi\cdot\nabla\theta \cos\theta 
\right)
\sin\theta\nabla\phi
-
\nabla
\left(|\nabla\mathbf{d}|^2\right)
\right]
\otimes \mathbf{d}
\right\|^2
\end{aligned}
\end{equation}
and
\begin{equation}
\begin{aligned}
&\|\nabla(\Delta\mathbf{d}+|\nabla\mathbf{d}|^2\mathbf{d})\|^2\\
=&
\left\|
\nabla
\left[
\left( \Delta\phi \sin\theta
+2\nabla\phi\cdot\nabla\theta \cos\theta 
\right)   \mathbf{d}^\perp_1
+(\Delta\theta- \sin\theta\cos\theta |\nabla\phi|^2 )\mathbf{d}^\perp_2
\right]
\right\|^2\\
=&
\left\|
\left[
\nabla\left( \Delta\phi \sin\theta
+2\nabla\phi\cdot\nabla\theta \cos\theta 
\right)
+
(\Delta\theta- \sin\theta\cos\theta |\nabla\phi|^2 )
\cos\theta \nabla\phi
\right]
\otimes  \mathbf{d}^\perp_1
\right\|^2 \\
+&
\left\|
\left[
\nabla(\Delta\theta- \sin\theta\cos\theta |\nabla\phi|^2 )
-\left( \Delta\phi \sin\theta
+2\nabla\phi\cdot\nabla\theta \cos\theta 
\right)\cos\theta\nabla\phi
\right]
\otimes \mathbf{d}^\perp_2
\right\|^2\\
+&
\left\|
\left[
(\Delta\theta- \sin\theta\cos\theta |\nabla\phi|^2 )\nabla\theta
+
\left( \Delta\phi \sin\theta
+2\nabla\phi\cdot\nabla\theta \cos\theta 
\right)
\sin\theta\nabla\phi
\right]
\otimes \mathbf{d}
\right\|^2.
\end{aligned}
\end{equation}

We just need to control three terms,
$\| |\nabla\mathbf{d}|^2 \sin\theta\nabla\phi\|$,
$\| |\nabla\mathbf{d}|^2 \nabla\theta  \|  $
and
$\|  \nabla \left(|\nabla\mathbf{d}|^2 \right) \|$,
by 
$\|\nabla(\Delta\mathbf{d}+|\nabla\mathbf{d}|^2\mathbf{d})\|$
and
$\| \Delta\mathbf{d}+|\nabla\mathbf{d}|^2\mathbf{d} \|$.
For $\| |\nabla\mathbf{d}|^2 \sin\theta\nabla\phi\|$
and
$\| |\nabla\mathbf{d}|^2 \nabla\theta  \|  $,
by Agmon inequality and Gagliardo-Nirenberg-Sobolev inequality \cite{ref6},
\begin{equation}
\| |\nabla\mathbf{d}|^2 \sin\theta\nabla\phi\|
+
\| |\nabla\mathbf{d}|^2 \nabla\theta  \|
\leq
2\| \nabla\mathbf{d}  \|^3_{L^6}
\leq
C\|\nabla\mathbf{d}\|
\left( \|\Delta\mathbf{d}\|+   \|\nabla\mathbf{d}\|  \right)^2
\leq
C\left( \|\Delta\mathbf{d}\|^2 +   1  \right).
\end{equation}
For $\|  \nabla \left(|\nabla\mathbf{d}|^2 \right) \|$,
we have
$$
\|  \nabla \left(|\nabla\mathbf{d}|^2 \right) \|
\leq
C\|\Delta\mathbf{d}\| \|\nabla\mathbf{d}\|_{L^{\infty}}
\leq
C\| \Delta\mathbf{d}  \| \|\nabla\Delta\mathbf{d}\|^{\frac{1}{2}}
\leq
\frac{1}{10} \|\nabla\Delta\mathbf{d}\|
+
C\| \Delta\mathbf{d}  \|^2
.
$$
We can get (2.17) from the above inequalities.
   $\Box$

\begin{remark}
For $k=1,2$,
we have the following identities. 
$$
\begin{aligned}
 \partial_k  \mathbf{d}
&= \left(\sin\theta \partial_k \phi \right)\mathbf{d}^\perp_1 
+
\left(\partial_k \theta    \right)  \mathbf{d}^\perp_2 \\
\partial_k  \mathbf{d}^\perp_1
&= -\left(\sin\theta \partial_k \phi \right)\mathbf{d} 
-\left(\cos\theta  \partial_k \phi\right)  \mathbf{d}^\perp_2   \\
\partial_k  \mathbf{d}^\perp_2
&= -  \left(\partial_k \theta \right) \mathbf{d}
+\left(\cos\theta \partial_k \phi\right) \mathbf{d}^\perp_1 
.
\end{aligned}
$$
\end{remark}

\begin{theorem}
Assuming that $\mathbf{d}:\mathbb{T}^3\rightarrow\mathbb{S}^2$ and
$\mathbf{d}\in C^{\infty}$
, 
we have the following estimate,
\begin{equation}
    \|\nabla^2\mathbf{d}\|^2
\leq
    C\left(
     \| \Delta\mathbf{d}+|\nabla\mathbf{d}|^2\mathbf{d} \|^3+1
     \right)   
\end{equation}
and
    \begin{equation}
    \|\nabla^3\mathbf{d}\|
\leq
    C\left(
    \| \Delta \mathbf{d} + |\nabla\mathbf{d}|^2 \mathbf{d} \|^\frac{9}{2}
    +
    \| \nabla(\Delta\mathbf{d}+|\nabla\mathbf{d}|^2\mathbf{d}) \|
    +1 \right),
\end{equation}
where $C$ is a constant depending on $\|\nabla\mathbf{d}\|$.
\end{theorem}

\paragraph{Proof.}
By Gagliardo-Nirenberg-Sobolev inequality \cite{ref6},
we have
$$
\begin{aligned}
\|\nabla\theta\|^4_{L^4}
&\leq
C \left(
\| \Delta\theta \|^3+1
\right) .
\end{aligned}
$$
Then we can obtain (2.21) 
by following a proof process similar to that of Theorem 2.1.

By Agmon inequality and Gagliardo-Nirenberg-Sobolev inequality \cite{ref6},
we have
$$
\begin{aligned}
\| \nabla\mathbf{d}  \|^3_{L^6}
&\leq
C
\left( \|\Delta\mathbf{d}\|^3 +1  \right),
\\
\|  \nabla |\nabla\mathbf{d}|^2  \|
&\leq
C\|\Delta\mathbf{d}\| \|\nabla\mathbf{d}\|_{L^{\infty}}\\
&\leq
C\|\Delta\mathbf{d}\|
\left(
\|\Delta\mathbf{d}\|^{\frac{1}{2}}
\|\nabla^3\mathbf{d}\|^{\frac{1}{2}}+1
\right)\\
&\leq
\frac{1}{10} \|\nabla\Delta\mathbf{d}\|
+
C\left(\| \Delta\mathbf{d}  \|^3 +1
\right)
.
\end{aligned}
$$
From (2.21), 
we have $\|\nabla^2\mathbf{d}\| \leq
    C\left(
     \| \Delta\mathbf{d}+|\nabla\mathbf{d}|^2\mathbf{d} \|^{\frac{3}{2}}+1
     \right)   $,
then we can obtain 
(2.22) 
by following a proof process similar to that of Theorem 2.2.
$\Box$

\begin{remark}
Because we are dealing with local coordinates on a sphere, 
instead of singularities,
we can address situations where angles on a closed curve lack single-valuedness.
It is helpful to get the higher-order estimates
and  the global well-posedness
of $(\mathbf{v},\mathbf{d})$ in the general Ericksen–Leslie system.
\end{remark}

{}	


\begin{thebibliography}{99}
	\itemsep=0pt
 
\bibitem{ref1}
ERICKSEN, J. L., Conservation laws for liquid crystals. Trans. Soc. Rheol. \textbf{5}, 23-34 (1961)
\bibitem{ref2}
HUANG, J., LIN, F., WANG, C., Regularity and Existence of Global Solutions to the Ericksen-Leslie System in $\mathbf{R}^2$. Commun. Math. Phys. \textbf{331}, 805-850 (2014)
\bibitem{ref3}
HUANG, X., A representation formula of incompressible liquid crystal flow and its applications. Nonlinearity \textbf{30}, 1911-1919 (2017)
\bibitem{ref4}
LEI, Z., LI, D., ZhANG, X., Remarks of global wellposedness of liquid crystal flows and heat flows of harmonic maps in two dimensions. Proc. Am. Math. Soc. \textbf{142}, 3801–3810 (2014)
\bibitem{ref5}
LESLIE, F. M., Some constitutive equations for liquid crystals. Arch. Ration. Mech. Anal. \textbf{28}, 265-283 (1968)
\bibitem{ref6}
NIRENBERG, L., On elliptic partial differential equations. Ann. Scuola Norm. Sup. Pisa Cl. Sci. \textbf{13}, 115–162 (1959)
\bibitem{ref7}
WANG, W., ZHANG, P., ZHANG, Z., Well-Posedness of the Ericksen-Leslie System. Arch. Ration. Mech. Anal. \textbf{210}, 837-855 (2013)

\end{thebibliography}
\end{document}